\documentclass[10pt]{amsart}
\usepackage{amssymb}
\usepackage{color}
\usepackage{graphicx}
\usepackage{float}
\usepackage[all,cmtip]{xy}
\usepackage{tikz}
\usetikzlibrary{matrix}
\usepackage{url}
\usepackage{subfig}
\usepackage{hyperref}
\mathchardef\mhyphen="2D

\newcommand{\std}{{\operatorname{std}}}

\newcommand{\Diff}{\operatorname{Diff}}

\newcommand{\SO}{\operatorname{SO}}
\newcommand{\U}{\operatorname{U}}

\newcommand{\Rot}{\operatorname{Rot}}
\newcommand{\tb}{\operatorname{tb}}

\newcommand{\Top}{\operatorname{Top}}

\newcommand{\Fr}{\operatorname{Fr}}

\newcommand{\R}{{\mathbb{R}}}
\newcommand{\C}{{\mathbb{C}}}
\newcommand{\Z}{{\mathbb{Z}}}

\newcommand{\NS}{{\mathbb{S}}}

\newcommand{\F}{{\mathcal{F}}}

\newcommand{\Leg}{{\mathfrak{Leg}}}

\newcommand{\Emb}{{\mathfrak{Emb}}}
\newcommand{\FEmb}{{\mathfrak{Emb}^{\Fr}}}

\newcommand{\wind}{\operatorname{winding}}
\newcommand{\writhe}{\operatorname{writhe}}

\newtheorem{theorem}{Theorem}

\newtheorem{proposition}[theorem]{Proposition}

\theoremstyle{definition}
\newtheorem{definition}[theorem]{Definition}

\newtheorem{remark}[theorem]{Remark}

\begin{document} 

\title{On the Legendrian realisation of loops of knots}

\subjclass[2020]{Primary: 53D10.}
\date{\today}

\keywords{}

\author{Javier Martínez-Aguinaga}

\begin{abstract}

We study the natural inclusion of the space of Legendrian embeddings in $(\NS^3,\xi_\std)$ into the space of smooth embeddings at the level of fundamental groups. T. Kálmán posed in \cite{kalman} the open question of whether for every fixed knot type $\mathcal{K}$ and Legendrian representative $\mathcal{L}$, the homomorphism $\pi_1(\mathcal{L})\to\pi_1(\mathcal{K})$ is surjective. We positively answer this question for infinitely many knot types $\mathcal{K}$ in the three main families (hyperbolic, torus and satellites) and every stabilised Legendrian representative in $(\NS^3,\xi_{std})$. 

\end{abstract}

\maketitle
\setcounter{tocdepth}{1} 

\section{Introduction}

The theory of Legendrian embeddings in the standard $3$-dimensional contact sphere $(\NS^3,\xi_\std)$ is a central topic in low-dimensional Contact Topology. The homotopy type of the space of Legendrian embeddings $\Leg(\NS^3,\xi_\std)$  in $(\NS^3,\xi_\std)$ is far from being completely understood, although much progress has been done. The homotopical properties of the natural inclusion of the space of Legendrian embeddings $\Leg(\NS^3,\xi_\std)$ into the space $\Emb(\NS^1,\NS^3)$ of smooth embeddings in $\NS^3$ provides information about how much the topology of these two spaces differ:

\begin{equation}\label{NaturalInclusion}
    i: \Leg(\NS^3,\xi_\std)\hookrightarrow \Emb(\NS^1,\NS^3).
\end{equation}

The induced map at the $\pi_0$-level has been extensively studied for many connected components of the space of Legendrian embeddings  \cite{Chekanov, EliashbergFraser, EtnyreHondaTorus, EtnyreHonda, FT, Ng}. It follows from the work of D. Fuchs  and S. Tabachnikov \cite{FT} that every smooth isotopy class of smooth embeddings admits a Legendrian representative. See also \cite[Thm 2.5]{Etnyre}.

As for the failure of injectivity at the $\pi_0$-level, note that this may happen purely for formal reasons since classical invariants (rotation number and Thurston-Bennequin invariant) can distinguish Legendrian isotopy classes representing the same smooth class. For a while, it was an open question whether Legendrian embeddings with the same knot type, rotation number and Thurston-Bennequin invariant were, in fact, Legendrian isotopic (see \cite[1.1.]{FT}).

Y. Chekanov \cite{Chekanov} introduced a Floer-Theoretical invariant that allowed him to find examples of Legendrian embeddings with the same formal invariants which were not Legendrian isotopic, albeit smoothly isotopic. This showed that the failure at the $\pi_0$-level was, in fact, more intricate that it may have been previously thought. Later on, more Legendrian classes were shown to satisfy this property by differents methods based on Convex surface Theory (see e.g.  \cite{EtnyreHonda}).

T. Kálmán \cite{kalman, kalman2} initiated the study of higher-homotopy groups of the space of Legendrian embeddings in the $3$-dimensional space. He constructed infinitely many loops of Legendrian positive torus embeddings that were contractible as loops of smooth embeddings but non-contractible in the space of Legendrian embeddings. As a consequence, he showed the lack of injectivity of the $\pi_1$-homomorphism induced by the natural inclusion (\ref{NaturalInclusion}) for infinitely many connected components.

\begin{theorem}[T. Kálmán, \cite{kalman}] There exist Legendrian knot types $\mathcal L$ so that for the corresponding smooth knot type $\mathcal K \supset \mathcal L$, the homomorphism $\pi_1(\mathcal L) \to \pi_1(\mathcal K)$ induced by the inclusion is not injective.
\end{theorem}

He then posed the open problem of whether this homomorphism, for fixed $\mathcal K$ and $\mathcal L$, is always surjective. We study this question in the present paper. Precisely, we show that this homomorphism is surjective for infinitely many knot types $\mathcal K$ in each of the three main families (hyperbolic, torus and satellites) and every stabilised representative. Furthermore, in the case of torus embeddings, not only we show that the result holds for infinitely many classes but we prove that it actually holds for all smooth torus knot types.

These results are stated in more generality throughout the article in Theorem \ref{mainTorus}, Theorem \ref{MainHyp} and Theorem \ref{mainPISat}. Let us encapsulate them in a unified (although less general) statement in the following Theorem \ref{mainPI1}.

\begin{theorem}\label{mainPI1}

The fundamental group homomorphism 
\begin{equation} \pi_1\left(\Leg(\NS^3,\xi_\std), \gamma\right)\rightarrow \pi_1\left(\Emb(\NS^1,\NS^3), \gamma\right)
\end{equation}

induced by the natural inclusion (\ref{NaturalInclusion})
is surjective for infinitely many smooth isotopy classes in each of the three main families (hyperbolic embeddings, $(p,q)$-torus embeddings and satellite embeddings) and every stabilised Legendrian representative $\gamma\in\Leg(\NS^3,\xi_\std)$. 

\end{theorem}

\begin{remark}
Note that we do not ask for the number of stabilisations to be large; if the considered Legendrian representative has only one stabilisation (either positive or negative), our result applies.
\end{remark}

\textbf{Acknowledgements:} I would like to thank Ryan Budney for teaching me about fundamental aspects of smooth embedding spaces during my visit to Victoria as well as in subsequent communications. I would also like to thank him for pointing out the existence of hyperbolic knots whose complements do not admit non-trivial hyperbolic isometries. Additionally, I would like to thank Eduardo Fernández and Francisco Presas for very interesting dicussions about Legendrian embedding spaces.

\section{Preliminaries on Legendrian embeddings}

\subsection{Legendrian embeddings in the standard $3$--dimensional space.}

Let us start introducing the standard contact structures both on $\R^3$ and on $\NS^3$.

The standard contact structure on $\R^3$  is given by $\xi_\std:=\ker(dz-ydx)$, whereas the standard contact structure on $\NS^3\subset\mathbb{C}^2$, which we also denote by $\xi_\std$ is given by 
\begin{equation}\label{ComplexTangencies}\xi_\std=T\NS^3\cap  JT\NS^3,\end{equation} where $J:T\mathbb{C}\to\mathbb C$ denotes the standard complex structure on $\C^2$.

\begin{remark}\label{contStereo}
    Both structures are closely related. In fact, $(\R^3,\xi_\std)$ is contactomorphic to $(\NS^3\backslash\{p\},\xi_{\std_{|\NS^3\backslash\{p\}}})$, for any choice of point $p\in\NS^3$ (\cite[Prop. 2.1.8]{GeigesBook} via the contact stereographic projection. 
\end{remark}

\begin{definition}
	A \textbf{Legendrian embedding} $\gamma:\NS^1\to(\NS^3,\xi_\std)$ is an embedding that is everywhere tangent to $\xi_\std$; i.e. $\gamma'(t)\in(\xi_\std)_{\gamma(t)}$ for every $t\in\NS^1$. Analogously, we say a $1$-dimensional embedded circle $L\subseteq (\NS^3,\xi_\std)$ is a \textbf{Legendrian knot} if it is the image of a Legendrian embedding.  
\end{definition}

Denote by $\widehat{\Leg}(M,\xi)$ the space of Legendrian knots of $(M,\xi)$ and by $\Leg(M,\xi)$ the space of Legendrian embeddings of $(M,\xi)$. The following relation holds $\widehat{\Leg}(M,\xi)=\Leg(M,\xi)/\Diff^+(\NS^1)$.

\subsubsection{Projections.} It is well known that we can recover a unique Legendrian curve in $(\R^3,\xi_\std)$ from curves in $\R^2$ that satisfy certain conditions. We will focus on the Lagrangian projection.

\begin{definition}
We define the \textbf{Lagrangian projection} as 
	\begin{center}
		$\begin{array}{rccl}
		\pi_{L}\colon & \R^3 & \longrightarrow & \R^2\\
		& (x,y,z)& \longmapsto & (x,y).
		\end{array}$
	\end{center}
\end{definition}
This projection has the property that maps immersed Legendrian curves in $(\R^3,\xi_\std)$ to immersed curves in $\R^2$. In addition, the $z$--coordinate can be recovered by integration:
\[z(t_1)=z(t_0)+\int_{t_0}^{t_1}y(s)x'(s)ds.\] 
Thus, in order for a closed curve in $\R^2$ to lift to a closed Legendrian, it is necessary that the closed immersed disk that bounds the curve has zero (signed) area.

\subsection{Classical invariants.}

There are three \textbf{classical invariants} of Legendrian embeddings that we will introduce, for simplicity, only in the context of $(\R^3,\xi_\std)$ and $(\NS^3,\xi_\std)$. The first one is the \textbf{smooth knot type} of the embedding, which is a purely topological invariant.

Let $(M,\xi)$ be $(\R^3,\xi_\std)$ or $(\NS^3,\xi_\std)$. Let $\gamma\in\Leg(M,\xi)$ be a Legendrian embedding.
\begin{definition} 
We define the \textbf{contact framing} of $\gamma$ as the trivialisation of its normal bundle $\nu(\gamma)$ given by the Reeb vector field along the knot.
\end{definition}

\begin{definition}
    We define the \textbf{topological framing} of $\gamma$ as the framing $\F_{\Top}$ of $\nu(\gamma)$ defined by a Seifert surface of $\gamma$.
\end{definition}

\begin{definition}\label{Def:tb}
The \textbf{Thurston-Bennequin invariant} of $\gamma$, denoted by $\tb(\gamma)$, is the twisting of the contact framing with respect to the topological framing.
\end{definition}

\begin{remark}
    The Thurston-Bennequin invariant does not depend on the choice of Seifert surface.
\end{remark}

Fix a global trivialisation of $\xi$. This election is unique up to homotopy since the space of maps from $M$ to $\NS^1$ is connected for the two particular manifolds that we are  studying. Thus, the derivative of the Legendrian embedding defines a map $\gamma':\NS^1\rightarrow\R^2\backslash\{0\}$.

\begin{definition}
We define the \textbf{rotation number} of $\gamma\in\Leg(M,\xi)$ as 
\[\Rot(\gamma)=\deg\gamma'.\]
\end{definition}

It follows that the rotation number is well defined and independent of the trivialisation of $\xi$. It is a well known result that the rotation number and the Thurston-Bennequin invariant of a Legendrian embedding $\gamma:\NS^1\to(\R^3,\xi_{std})$ can be easily described in terms of its Lagrangian projection $\pi_L(\gamma)$.

\begin{proposition}
    The rotation number of $\gamma$ can be computed as the winding number of its Lagrangian projection; i.e. $\Rot(\gamma)=\wind(\pi_L\circ\gamma)$.
\end{proposition}

\begin{proposition}
    The Thurston-Bennequin of $\gamma$ coincides with the writhe of its Lagrangian projection; i.e. $\tb(\gamma)=\writhe(\pi_L\circ\gamma)$.
\end{proposition}

\begin{remark}\label{reeb}
Note that the Reeb vector field $\partial_z$ is transverse to the Lagrangian projection.
\end{remark}

\section{Some homotopical relation between smooth and Legendrian embeddings}

We will show some relevant relations between the space of Legendrians and the space of long Legendrians focusing on fundamental groups. Let us introduce the situation in the smooth setting first.

\subsection{Smooth embeddings in $\NS^3$ and long embeddings in $\R^3$.}

\begin{definition}
    A \textbf{long embedding}  is an embedding $\gamma:\mathbb R\to \mathbb R^3$ such that $\gamma(t)=(t,0,0)$ for $|t|>1$. We denote the space of long embeddings by $L\Emb(\R,\R^3)$. 
\end{definition}

\begin{remark}
    Note that, for $|t|>1$, all long embeddings $\gamma(t)$ coincide  along the axis $\{(t,0,0):t\in\R\}\subset\R^3$, which is often called the ``long axis'' in the literature.
\end{remark}

This space is homotopically equivalent to the following space.

\begin{definition}\label{LongEmbeddings}
    Write $\Emb_{N,jN}(\NS^1,\NS^3)$ for the space of smooth embeddings $\gamma(t)$ in $\NS^3$ passing at $t=0$ through the north pole $N=(1,0)\in\C^2$ with fixed derivative $j\cdot N=(0,1)\in\C^2$ up to scaling (where we are using quaternionic notation). Precisely:
\[\Emb_{N,jN}(\NS^1,\NS^3)=\left\{\gamma\in\Emb(\NS^1,\NS^3):\gamma(0)=N,\ \frac{\gamma'(0)}{||\gamma'(0)||}=j\cdot N\right\}.\] 
\end{definition}

\begin{remark}\label{NoAmbiguity}
The homotopical identification of the spaces $L\Emb(\R,\R^3)$ and  $\Emb_{N,jN}(\NS^1,\NS^3)$ is standard in the literature of embedding spaces. See, for instance, See Proposition 4.1 and its proof in \cite{BudCohen}. Alternately, the reader can consult \cite[Thm 2.1]{Budney2} or the second part of the proof of \cite[Corollary 6]{Budney2}. This identification can be geometrically understood as follows. The one-point compactification of  $\R^3$ into $\NS^3$ allows to ``close'' any long embedding into a an embedding in $\NS^3$ with fixed evaluation at $0\in\NS^1$ and fixed derivative up to scaling. Likewise, we can revert the process by taking the stereographic projection \cite[Prop. 4.1]{BudCohen}, so that $L\Emb(\R,\R^3)$ and $\Emb_{N,jN}(\NS^1,\NS^3)$ are naturally identified.

By means of this identification, we will refer to both spaces indistinguishably as the space of long embeddings. Likewise, we will often work with both spaces, passing from one to the other at our convenience. 
\end{remark}

\begin{remark}[Notation]\label{notation}
Henceforth, when we write a superscript $\gamma$ over any embedding space in this article, we are denoting the connected component of such space containing 
 the embedding $\gamma$. For instance, $\Emb^\gamma(\NS^1,\NS^3)$ denotes the connected component of the space of smooth embeddings $\Emb(\NS^1,\NS^3)$ containing the embedding $\gamma$.
\end{remark}

Following \cite{Budney2, Budney3, BudCohen}, we have the following well defined smooth fibration
\begin{equation}\label{eq:LongEmbeddingsFibration}
\Emb_{N,jN}(\NS^1,\NS^3)\hookrightarrow\Emb(\NS^1,\NS^3)\rightarrow V_{4,2}=\SO(4)/\SO(2)
\end{equation}
where the projection onto the base $V_{4,2}$ is just the $1$-jet projection with normalised derivative; i.e.  $\rho:\Emb(\NS^1,\NS^3)\to V_{4,2}$, $\gamma\mapsto(\gamma(0)|\frac{\gamma'(0)}{||\gamma'(0)||})$. Note that the fiber over $(N,jN)\in V_{4,2}$ is precisely $\Emb_{N,jN}(\NS^1,\NS^3)$. Moreover, the fiber $\Emb_{N,jN}(\NS^1,\NS^3)$ embeds into $\Emb(\NS^1,\NS^3)$ just by the natural inclusion. This fibration relates the topology of the space of smooth embeddings in $\NS^3$ with the topology of the space of long embeddings-

 Fibration (\ref{eq:LongEmbeddingsFibration}) induces a long exact sequence at the level of homotopy groups. By focusing at the $\pi_1$-level, we can write:
\begin{equation}\label{eq:LongEmbeddingsFibrationpi1}
\pi_1\left(\Emb_{N,jN}(\NS^1,\NS^3)\right)\rightarrow\pi_1\left(\Emb(\NS^1,\NS^3)\right)\rightarrow \pi_1(V_{4,2})
\end{equation}

\subsection{Legendrian embeddings and long legendrians}
We have an analogous situation in the contact setting, as already pointed out in \cite{FMP3} and also treated in \cite{FM}.

The homotopy type of the space of Legendrian embeddings in the standard $(\NS^3,\xi_\std)$ is intimately related to the space of long Legendrian embeddings $\Leg_{N,jN}(\NS^3,\xi_\std):=\Emb_{N,jN}(\NS^1,\NS^3)\cap\Leg(\NS^3,\xi_\std).$ Analogously, we write $L\Leg(\R,\R^3)$ for the subspace of $L\Emb(\R,\R^3)$ satisfying the Legendrian condition. 
When we focus on fundamental groups, we have the following statement.

\begin{proposition}\label{PropSurj1} The group homomorphism induced by the natural inclusion at the $\pi_1$-level 
\begin{equation}\label{Embpi1} \pi_1\left(\Leg(\NS^3,\xi_\std), \gamma\right)\rightarrow \pi_1\left(\Emb(\NS^1,\NS^3), \gamma\right)
\end{equation}
is surjective if the associated morphism at the level of long components is surjective.\begin{equation}\label{Legp1}
\pi_1\left(\Leg_{N,jN}(\NS^3,\xi_\std), \gamma\right)\rightarrow\pi_1\left(\Emb_{N,jN}(\NS^1,\NS^3), \gamma\right).
\end{equation} 
\end{proposition}

\begin{proof}
Since $\pi_1(V_{4,2})=\{0\}$, then by the long exact sequence (\ref{eq:LongEmbeddingsFibrationpi1}), we get that the group homomorphism $\pi_1\left(\Emb_{N,jN}(\NS^1,\NS^3), \gamma\right)\rightarrow \pi_1\left(\Emb(\NS^1,\NS^3), \gamma\right)$ is surjective. Therefore, since $\Leg_{N,jN}(\NS^3,\xi_\std)$ is a subspace of $\Leg(\NS^3,\xi_\std)$, then the surjectivity of the group homomorphism (\ref{Legp1}) readily implies the surjectivity of the surjectivity of the group homomorphism (\ref{Embpi1}).
\end{proof}

\begin{remark}
    Be aware that we are identifying the embedding $\gamma \in\Leg(\NS^3,\xi_\std)$ with its corresponding long embedding $\Leg_{N,jN}(\NS^3,\xi_\std)$. This is implicit in the used notation.
\end{remark}

\section{Two canonical loops in the space of embeddings}

This Section is devoted to the study of two classical loops in the space of smooth embeddings; namely, Gramain's loop and Fox-Hatcher loop. As we will see, these two loops are canonically well defined at the homotopy level for every connected component of the space of smooth embeddings. We will study how to realise these two loops in the Legendrian setting.

\subsection{Gramain's loop}

Gramain's loop, first studied by A. Gramain in \cite{Gramain}, yields a non-trivial element in the fundamental group of the space of long embeddings $\Emb^\gamma_{N,jN}(\NS^1,\NS^3)$ for all non-trivial embeddings $\gamma$. It plays a crucial role in the homotopical description of embedding spaces and we will show how to realise it by stabilised Legendrian embeddings.

\begin{definition}
	The \textbf{Gramain loop} associated to a long embedding $\gamma\in\Emb_{N,jN}(\NS^1,\NS^3)$ is the loop 
	$$ gr(\gamma)=B_\theta\cdot\gamma, \quad \theta\in[0,2\pi],\quad\text{ where }$$
 
	\[B_\theta := \left( {\begin{array}{cccc}
		1 & 0 & 0 & 0\\
		0 & \cos \theta & 0 &-\sin\theta \\
		0 & 0 & 1  & 0 \\
		0 & \sin\theta & 0 & \cos\theta
		\end{array} } \right)\in \SO(4).\]	
\end{definition}

Equivalently, it can be defined as a full rotation along the long axis of the long embedding $\gamma\in L\Emb(\R,\R^3)$ under the identification (\ref{NoAmbiguity}). Note that this matrix fixed both $N=(1,0)$ and $JN=(0,1)$ or alternatively, using real coordinates, $(1,0,0,0)$ and $(0,0,1,0)$, respectively.

\subsubsection{Legendrian realisation of Gramain's loop for stabilised Legendrians} 

The following Legendrian realisation of Gramain's loop was inspired to the author by a geometric depiction of Gramain's smooth loop in the book \cite[2.3]{Fiedler} of T. Fiedler.

\begin{proposition}\label{Prop:GrRealisation}
    Let $\gamma\in \Leg_{N,jN}(\NS^3, \xi_\std)$ be a long Legendrian embedding with at least one stabilisation. Then, Gramain's loop and its inverse can be realised as loops $\gamma^{\theta}$ of long Legendrian embeddings based at $\gamma$.
\end{proposition}

\begin{proof}
    
Consider $\gamma:\NS^1\to\NS^3$ a long Legendrian embedding with a positive stabilisation (the argument will be symmetric for a negative stabilisation). Equivalently, $\gamma$ can be regarded as a connected sum $\gamma:=\tilde{\gamma}\#\beta$ where $\tilde{\gamma}$ is some long legendrian embedding and $\beta$ denotes a long legendrian unknot with one positive stabilisation. We can diagrammatically represent $\gamma$ as in Figure \ref{FIG1}, where the first box $K_{\tilde{\gamma}}$ depicts the long embedding $\tilde{\gamma}$ and the second one, $K_{\beta}$, encodes the positive stabilisation.

\begin{figure}[h!]
	\centering
	\includegraphics[width=0.6\textwidth]{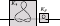}
	\caption{Every stabilised long legendrian embedding $\gamma$ can be regarded, up to Legendrian isotopy, as a connected sum $\gamma:=\tilde{\gamma}\#\beta$ where $\tilde{\gamma}$ is some long legendrian embedding and $\beta$ denotes a long legendrian unknot with one positive stabilisation.\label{FIG1}}
\end{figure}

Consider then the following nice realisation of Gramain's loop for $\gamma$. We define it locally on a local Darboux chart disjoint from $N$. We work on a Lagrangian projection, where the following loop is defined on a small arc not containing the origin $\gamma(0)$, whereas everything else remains constant away from this arc.

On the first step, after adjusting sizes of the boxes appropriately, we pass  the box $K_{\tilde{\gamma}}$ all through $K_\beta$ until both boxes have exchanged positions. This is depicted in Figure \ref{FIG2}. On a second step, we pass the box $K_\beta$ all through $K_{\tilde{\gamma}}$ until both boxes come back to their initial positions.

\begin{figure}[h]
	\centering
	\includegraphics[width=0.85\textwidth]{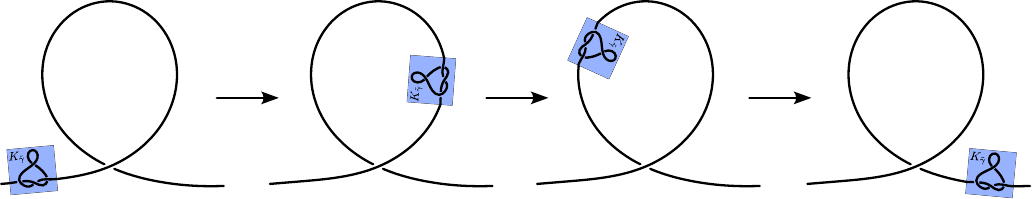}
	\caption{First part of the loop where the $K_{\tilde{\gamma}}-$box passes all through the stabilisation $\beta$. Since $\tb(\tilde\gamma\#\beta)=\tb(\tilde\gamma)-1$ (and, thus, the stabilisation adds a full rotation of the contact framing with respect to the topological framing), this first part of the loop corresponds to a full $2\pi$-rotation of the $K_{\tilde{\gamma}}-$box around the long axis. When concatenated with the second part of the loop, this yields a Legendrian realisation of Gramain's loop. \label{FIG2}}
\end{figure}

Since $\tb(\tilde\gamma\#\beta)=\tb(\tilde\gamma)-1$, and the Thurston-Bennequin invariant measures the twisting of the contact framing with respect to the topological framing (Definition \ref{Def:tb}), it readily follows that the first part of this loop is smoothly isotopic to performing a full rotation of the box $K_{\tilde{\gamma}}$ along the long axis and, thus, realises Gramain's loop for $\gamma$. The second part (where the kink passes all trough the long knot) is smoothly trivial since it consists of taking a stabilisation and passing it all along the long knot. Clearly this can be undone by parametrically killing this kink with a Reidemeister I move.

Finally, note that in case we had started with a negative stabilisation instead of a positive one, this same process would have yielded the realisation of the inverse of Gramain's loop. Therefore, the realisation of its inverse yields Gramain's loop in this other case as well. Finally,  note that in order to realise the respective inverses in both cases we just need to perform the loop in the reverse order. \end{proof}

\begin{remark}
    The loop described in Proposition \ref{Prop:GrRealisation} is a particular instance of the ``Legendrian pulling one knot through another'' family of loops of Legendrian embeddings introduced in \cite{FMP4}.
\end{remark}

\subsection{The Hatcher-Fox map}

We will now study the Hatcher-Fox loop, which is well defined for any given embedding together with a choice of framing. In particular, we will work with framed embeddings along the Subsection.

\begin{definition}
    Denote by $\FEmb(\NS^1,\NS^3)$ the space of framed embeddings of the circle $\NS^1$ into $\NS^3$, i.e. the space of pairs $(\gamma,\tau)$ where $\gamma\in\Emb(\NS^1,\NS^3)$ and $\tau=(\tau_1,\tau_2)$ is a trivialisation of the normal bundle of $\gamma$.
    
\end{definition}

\begin{remark}There is a canonical choice (up to homotopy) of trivialisation of the normal bundle of an embedding $\gamma$ determined by a Seifert surface $\Sigma$ of the embedding. We call this trivialisation the \textbf{topological framing}, and we denote it by $\tau_{Top}$.
\end{remark}
Given an embedding $\gamma\in\Emb(\NS^1,\NS^3)$ denote by $\Fr(\gamma)$ the space of trivialisations of the normal bundle of $\gamma$. If we trivialise the normal bundle of $\gamma$ by means of $\tau_{Top}$ we can identity $\Fr(\gamma)=L\NS^1$. Thus, $\pi_0(\Fr(\gamma))=\pi_0(L\NS^1)=\pi_1(\NS^1)=\Z$.

Given a framed embedding $(\gamma,\tau)$ we can build a loop of matrices in $\SO(4)$ defined by 

\[ A_{\gamma,\tau}(\theta) = \left( \gamma(\theta) \, \middle| \, \tau_1(\theta) \, \middle| \, \frac{\gamma'(\theta)}{||\gamma'(\theta)||} \, \middle| \, \tau_2(\theta) \right),\quad \theta\in\NS^1\simeq [0,2\pi].\]

Observe that $\tilde{\gamma}^\theta(t)=A_{\gamma,\tau}(\theta)^{-1}\cdot\gamma(t+\theta)$ satisfies that 

\begin{itemize}
	\item[i)] $\tilde{\gamma}^\theta(0)=A_{\gamma,\tau}(\theta)^{-1}\cdot\gamma(\theta)=(1,0,0,0)=N$
	\item[ii)] $(\tilde{\gamma}^{\theta})'(0)=A_{\gamma,\tau}^{-1}\cdot\gamma'(\theta)=||\gamma'(\theta)||\cdot(0,0,1,0)=||\gamma'(\theta)||\cdot j\cdot N.$
\end{itemize}

\begin{remark}\label{RmkFH}
We have, thus, that the embedding $\gamma^\theta(t)=\tilde{\gamma}^\theta\left(\frac{1}{||\gamma'(\theta)||}\cdot t\right)$ is, in fact, a long embedding, i.e. an element of the space $\Emb_{N,jN}(\NS^1,\NS^3)=\lbrace\beta\in\Emb(\NS^1,\NS^3):\beta(0)=N,\beta'(0)=j\cdot N\rbrace$.
\end{remark}

This construction is canonical and defines a continuous map 

\begin{equation}\label{HFMap}\mathcal{HF}:\FEmb(\NS^1,\NS^3)\rightarrow \pi_1\left(\Emb_{N,jN}(\NS^1,\NS^3)\right),\quad (\gamma,\tau)\mapsto\gamma^\theta,
\end{equation}

that we call the \textbf{Hatcher-Fox Map}.

\begin{definition}\label{SmoothHF}
	Given an embedding $\gamma\in\Emb(\NS^1,\NS^3)$ a \textbf{Hatcher-Fox loop} associated to $\gamma$ is any loop of the form
	$$ hf(\gamma)=\mathcal{HF}(\gamma,\tau_{Top}).$$
 
\end{definition}

\begin{remark}
Observe that the homotopy class of a Hatcher-Fox loop is well defined.
\end{remark}

\begin{proposition}[Hatcher]\label{GramainTorus}
	If $\gamma$ is a parametrised torus knot, then $hf(\gamma)=gr(\gamma)$. 
\end{proposition}

\begin{remark}
    
Consider a framed embedding $(\gamma,\tau)$. We assume without loss of generality that $\gamma\in\Emb_{N,jN}(\NS^1,\NS^3)$ and that $\tau(0)=(e_2,e_4)$. Use a topological framing of $\gamma$ to identify $\Fr(\gamma)=L\NS^1$. Then, $\tau=k\in\pi_0(L\NS^1)=\mathbb{Z}$. It follows that 

\begin{equation}\label{AddTwist}\mathcal{HF}(\gamma,\tau)=hf(\gamma)+k \cdot gr(\gamma)\in\pi_1(\Emb_{N,jN}(\NS^1,\NS^3),\gamma).
\end{equation}

\end{remark}

A. Hatcher provides a geometric depiction of the Hatcher-Fox loop in \cite{hatcher1}. He starts with a knot in the $3$-dimensional space and places a small bead on the knot (which acts as a wire throughout which the bead can move). By assumming that the endpoints of the arc inside the bead are  antipodal, the complement of the knot outside the bead can be regarded as a long embedding. 

Then by moving the bead all around the knot until it comes back to the original starting point, we get a loop of long embeddings. As Equation (\ref{AddTwist}) shows, changing the frame of $\gamma$ (as a framed embedding) is equivalent to adding or substracting twists of such  bead around the wire while performing the loop; i.e. those twists correspond to adding Gramain's loops or its inverses, respectively.

Note that, as stated in \cite{HK}, we can get the following equivalent depiction of Hatcher-Fox loop. Take a long knot $\gamma\in L\Emb(\R,\R^3)$ and consider its one-point compactification $\hat\gamma\in\Emb(\NS^1,\NS^3)$. Then one can move $\infty$ along the trace of $\hat\gamma$ and stereographically project to $\R^3$ accordingly with the frame as in Remark \ref{RmkFH}. See Figure \ref{FIG1} for such a representation in the Legendrian setting.

Let $\gamma\in\Leg(\NS^3,\xi_\std)$ be a Legendrian embedding. Then $\gamma$ is naturally framed by means of the \textbf{contact framing} 
\begin{equation}\label{ContactFraming}
\tau_{Cont}:=(i\gamma,i\gamma').
\end{equation} 

This defines an inclusion $j_{Cont}:\Leg(\NS^3,\xi_\std)\hookrightarrow\FEmb(\NS^1,\NS^3),\quad\gamma\mapsto(\gamma,\tau_{Cont})$.

It follows also that $A_{\gamma,\tau_{Cont}}\in\U(2)$ and, thus, preserves the contact structure $\xi_\std$ and, in particular, $\mathcal{HF}(\gamma,\tau_{Cont})\in  \Leg_{N,jN}(\NS^3,\xi_\std)$ is a loop of long Legendrian embeddings. 

\begin{definition}
	We define the \textbf{Legendrian Hatcher-Fox Map} as the composition 
\[\mathcal{HF}_{\Leg}=\mathcal{HF}\circ i_{Cont}:\Leg(\NS^3,\xi_\std)\rightarrow \Leg_{N,jN}(\NS^3,\xi_\std).\]
 
	\begin{definition}\label{LegHF}
	 The \textbf{Legendrian Hatcher-Fox loop} of $\gamma$ is defined as the loop 
 
 \[hf_{L}(\gamma)=\mathcal{HF}_{\Leg}(\gamma).\]
 
	\end{definition}
\end{definition}

\begin{remark}
By definition, $\tau_{Cont}=\tb(\gamma)\in\pi_0(\Fr(\gamma))=\mathbb{Z}$. In particular, the following relation holds:
\begin{equation}\label{relationTB} hf_{L}(\gamma)=hf(\gamma)+\tb(\gamma) \cdot gr(\gamma)\in\pi_1(\Emb_{N,jN}(\NS^1,\NS^3)).
\end{equation}
\end{remark}
Let us encapsulate the main constructions of this Subsection in the following result.

\begin{theorem}\label{Grealisations}
Let $\gamma\in\Leg_{N,jN}(\NS^3,\xi_\std)$ be a long Legendrian embedding with at least one stabilisation (either positive or negative). Then both Gramain's loop $[\mathcal G]\in\pi_1\left(\Emb_{N,jN}(\NS^1,\NS^3)\right)$ and Hatcher-Fox loop  $[\mathcal{HF}]\in\pi_1\left(\Emb_{N,jN}(\NS^1,\NS^3)\right)$ admit Legendrian representatives in the component of $\gamma$. More precisely, they possess preimages by the natural homomorphism \begin{equation*}
\pi_1\left(\Leg_{N,jN}(\NS^3,\xi_\std), \gamma\right)\rightarrow\pi_1\left(\Emb_{N,jN}(\NS^1,\NS^3),\gamma\right).
\end{equation*} 
\end{theorem}

\begin{proof}
    The realisation of Gramain's loop follows from Proposition \ref{Prop:GrRealisation}. On the other hand, we showed that the Legendrian Hatcher-Fox loop can be realised on every connected component of the space of Legendrian embeddings, without the need of restricting to stabilised Legendrian classes.
    
    Recall that the Legendrian Hatcher-Fox loop (Definition \ref{LegHF}) differs with the (topological) Hatcher-Fox loop (Definition \ref{LegHF}) precisely on a total number of $\tb(\gamma)$ Gramain's loops. Since $\gamma$ has at least one stabilisation, we can ``cancel'' those additional $\tb(\gamma)$-Gramain's loops by realising their inverses (Proposition \ref{Prop:GrRealisation}) and thus get the desired realisation for the Hatcher-Fox loop.
\end{proof}

\subsection{Loops of $(p,q)$-torus embeddings}

Smooth long $(p,q)$-torus embeddings present the particularity that their Hatcher-Fox loop coincides (at the homotopy level) with their Gramain's loop \cite{hatcher3}. Let us focus on this particular type of embeddings along the present Subsection.

\begin{remark} If $\gamma$ is a Legendrian torus $(p,q)$--embedding, then, by Proposition \ref{GramainTorus},

$$ hf_L (\gamma)=(1+\tb(\gamma))\cdot gr(\gamma)\in\pi_1(\Emb_{N,jN}(\NS^1,\NS^3)).$$

\begin{figure}[h]
	\centering
	\includegraphics[width=1\textwidth]{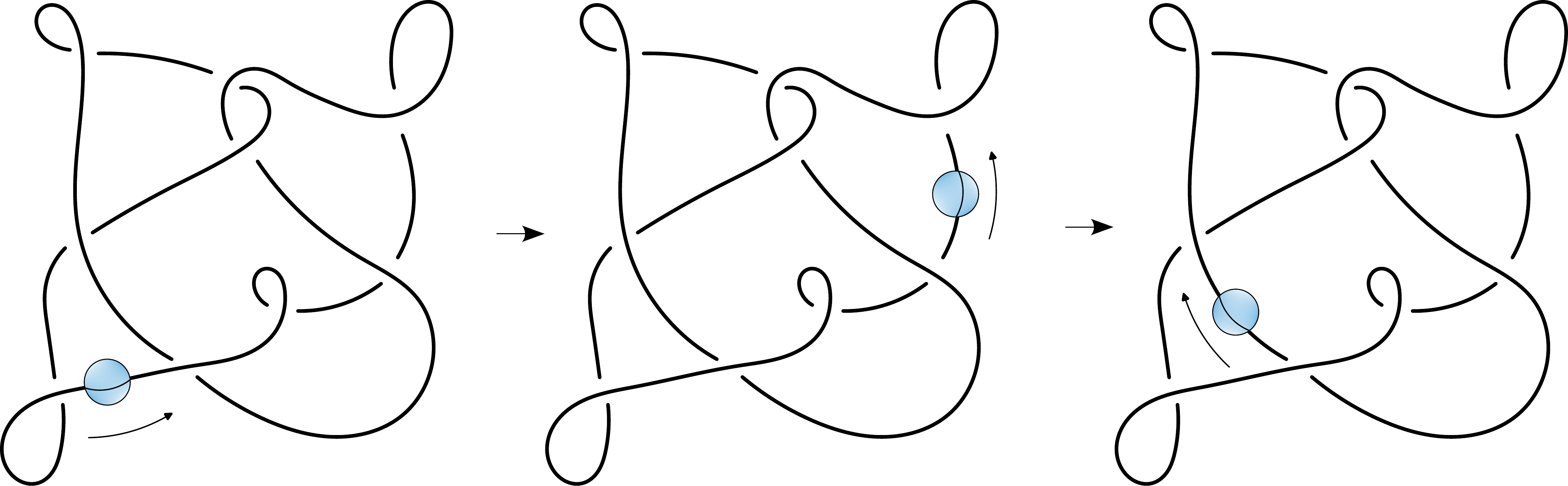}
	\caption{Legendrian realisation of Legendrian Hatcher-Fox loop  $hf(\gamma)+\tb(\gamma) \cdot gr(\gamma)$ for one of Chekanov's $m_{5_2}$  Legendrian knots $\gamma$. Geometrically, this loop can be thought as the loop of long embeddings arising from reparametrising $\gamma$ while simulatenously applying the Hatcher-Fox map defined in (\ref{HFMap}) by means of the contact framing (\ref{ContactFraming}). The blue pearl can be identified with an arbitrarily small neighborhood of $\infty\in(\NS^3,\xi_\std)$, the point from which we can stereographically project everything onto $(\R^3,\xi_\std)$ and thus obtain the aforementioned loop of long legendrian embeddings. \label{FIG1}}
\end{figure}

If $\gamma$ is a positive torus knot with maximal Thurston--Bennequin invariant $\tb(\gamma)=pq-p-q$, the previous formula takes the following explicit form:

\begin{equation} hf_L (\gamma)=(1+pq-p-q)\cdot gr(\gamma)=(p-1)(q-1)\cdot gr(\gamma)\in\pi_1(\Emb_{N,jN}(\NS^1,\NS^3)).\end{equation}

\end{remark}	

Note that, since we already showed how to realise Gramain's loop (as well as its inverse loop) by stabilised Legendrians, we can readily conclude the statement concerning torus embeddings in Theorem \ref{mainPI1} from the Introduction.

\begin{theorem}\label{mainTorus}

The fundamental group homomorphism 
\begin{equation} \pi_1\left(\Leg(\NS^3,\xi_\std), \gamma_{p,q}\right)\rightarrow \pi_1\left(\Emb(\NS^1,\NS^3), \gamma_{p,q}\right)
\end{equation}

induced by the natural inclusion (\ref{NaturalInclusion})
is surjective for every stabilised Legendrian $(p,q)$-torus embedding $\gamma_{p,q}\in\Leg(\NS^3,\xi_\std)$.
\end{theorem}

\begin{proof}
     By Proposition \ref{PropSurj1}, it suffices to prove the analogous statement at the level of long-components. Since the fundamental group of the space of long torus embeddings is generated by Gramain's loop \cite{hatcher2}, then the result just follows by Proposition \ref{Prop:GrRealisation}.
\end{proof}

\begin{remark}
    Note that the proper subspace of $\Leg(\NS^3,\xi_\std)$ consisting of positive Legendrian torus embeddings with maximal Thurston-Bennequin invariant has been studied by E. Fernández and H. Min in \cite{FM}. They have determined its homotopy type and checked that previous loops of Legendrian torus embeddings described in \cite{kalman} already generated the long component at the $\pi_1$-level. Our results extends the $\pi_1$-surjectivity statement for the whole class of stabilised embeddings, including all negative torus embeddings.
\end{remark}

\section{The subgroup $\mathfrak{G}$ and hyperbolic isometries}\label{symmetries}

A subgroup of interest in the study of embedding spaces is the subgroup of $\mathfrak G < \pi_1\left(L\Emb(\R,\R^3)\right)$ generated by Gramain's loop and Hatcher-Fox loop.

\begin{definition}
    We define the group $\mathfrak G < \pi_1\left(L\Emb(\R,\R^3)\right)$ as the subgroup generated by Gramain's loop and Hatcher-Fox loop. We will often denote $\mathfrak G_\gamma$ to refer to that same group when we want to make explicit the specific choice of basepoint $\gamma\in L\Emb(\R,\R^3)$.
\end{definition}

\begin{remark}
    Since these two loops are canonically defined up to homotopy for every long connected component of the space of long embeddings, this subgroup is then well defined for any connected component of the space $L\Emb(\R,\R^3)$.
\end{remark}

\begin{remark}
    This Definition readily extends to the space $\Emb_{N,jN}(\NS^1,\NS^3)$ by the obvious identification. Recall that we will indistinctively refer to both spaces $L\Emb(\R,\R^3)$ and $\Emb_{N,jN}(\NS^1,\NS^3)$  as the space of long embeddings and pass from one to another without any loss of generality nor ambiguity, as we already stated in Remark \ref{NoAmbiguity}.

\begin{remark}
    Theorem \ref{Grealisations} implies that the fundamental group of any stabilised component $L\Leg^\gamma(\R,\R^3)$ (respectively, $\Leg^\gamma_{N,jN}(\NS^3,\xi_\std)$) contains a copy of $\mathfrak{G}_\gamma$; i.e.
    
    \[
    \mathfrak{G}_\gamma < \pi_1\left(\Leg^\gamma_{N,jN}(\NS^3,\xi_\std) \right).
    \]
\end{remark}

    We will introduce the following group of transformations associated to a hyperbolic long embedding by regarding it, for convenience, as an embedding in $\NS^3$. We follow \cite{Budney3}. 
\end{remark}

\begin{definition}
Define $B_\gamma$ as the
group of hyperbolic isometries of the complement of a smooth long embedding $\gamma\in\Emb_{N,jN}(\NS^1,\NS^3)$ 
which:
\begin{itemize}
    \item[i)] extend to orientation-preserving
diffeomorphisms of $\NS^3$,
\item[ii)] preserve the orientation of $\NS^3$ and
\item[iii)] preserve $\gamma$ (as an embedded oriented curve) and its orientation.
\end{itemize}
\end{definition}

\begin{remark}\label{R1}
It follows from the work \cite{Budney3} of Ryan Budney that the group $\mathfrak G$ has finite index within  $\Emb^\gamma_{N,jN}(\NS^1,\NS^3)$ when $\gamma$ is hyperbolic and, furthermore, the index of $\mathfrak G$ coincides with the order of the group $B_\gamma$. 
\end{remark}

\begin{remark}\label{R2} There are infinitely many hyperbolic embeddings in $\NS^3$ that have no symmetries (i.e. there are not orientation preserving diffeomorphisms $f:\NS^3\to\NS^3$ of finite order that leaves the image of $\gamma$ invariant). The reader can consult, for instance, the work of M. Teragaito \cite[Thm. 1.2 \& Prop. 5.4]{Teragaito} for explicit examples. Also, as pointed out by Teragaito in \cite{Teragaito}, earlier work by Deruelle, Miyazaki, and Motegi \cite{DMM} also provides infinitely many explicit examples of hyperbolic knots without symmetries. These are constructed from the $(-3,3,5)$-pretzel knots. See \cite{DMM} for further details.

Therefore, all the aforementioned examples also provide infinitely many instances of embeddings $\gamma\in \Emb_{N,jN}(\NS^1,\NS^3)$ for which the group $B_\gamma$ is trivial.
\end{remark}

Putting together Remark \ref{R1} and Remark \ref{R2}. we readily get the following Proposition \ref{FGPI1}, which states that there are infinitely many hyperbolic long embeddings for which the associated Hatcher-Hatcher and Gramain's loops generate their whole fundamental group.

\begin{proposition}\label{FGPI1} 
    There are infinitely many hyperbolic long embeddings $\gamma\in L\Emb(\R,\R^3)$ satisfying that the natural group inclusion homomorphism
   \[i: \mathfrak{G}_\gamma\hookrightarrow
   \pi_1\left( L\Emb(\R,\R^3),\gamma\right)
   \]
    is a group isomorphism and, thus, both groups coincide. 
\end{proposition}

Finally, let us state the main Theorem of the Subsection (which corresponds to the statement concerning Hyperbolic embeddings in Theorem \ref{mainPI1} from the Introduction). 

\begin{theorem}\label{MainHyp}
The fundamental group homomorphism 
\begin{equation} \pi_1\left(\Leg(\NS^3,\xi_\std), \gamma\right)\rightarrow \pi_1\left(\Emb(\NS^1,\NS^3), \gamma\right)
\end{equation}

induced by the natural inclusion (\ref{NaturalInclusion})
is surjective for infinitely many smooth  hyperbolic embeddings and every stabilised Legendrian representative $\gamma\in\Leg(\NS^3,\xi_\std)$.
\end{theorem}

\begin{proof}
    By Proposition \ref{PropSurj1}, it suffices to prove the analogous statement at the level of long-components. By restricting ourselves to the infinitely many smooth hyperbolic components considered in Proposition \ref{FGPI1}, it suffices to generate the subgroup $\mathfrak{G}_\gamma$. Therefore, the result follows by Theorem \ref{Grealisations}.
\end{proof}

\section{Legendrian connected-sums and associated loops}

\subsection{The smooth case}

R. Budey provided in \cite{Budney3}  a detailed description of the whole homotopy type of the space $L\Emb(\R,\R^3)$ of long embeddings in $\mathbb{R}^3$ or, equivalently, of the space $\Emb_{N,jN}(\NS^1,\NS^3)$. This result is recursive; i.e. it describes the global homotopy  of all the connected  components of the space (precisely, the fundamental group since each component is a $K(\pi, 1)$, \cite{hatcher1, hatcher3}) iteratively building on simpler components. 

In particular, this includes a structural theorem for the connected sum of any finite number of prime components (which was first proven in \cite{Budney1} by means of the action of little cubes in the space of long embeddings). Let us provide a precise statement.

\begin{theorem}[R. Budney, \cite{Budney1, Budney3}]
    Let $\gamma$ be the connected sum of  $n\geq 2$ prime long embeddings $\gamma_1,\cdots, \gamma_n$.
    Then  $L\Emb^\gamma(\R,\R^3)\simeq C_2(n)\times_{\Sigma_\gamma}\prod_{i=1}^{n}  L\Emb^{\gamma_i}(\R,\R^3)$.
\end{theorem}

\begin{remark}[Notation]
We have used (almost) the same notation as \cite{Budney3}. 
    Here $C_2(n)$ denotes the configuration space of $n$ points 
    
    \[C_2(n) := \lbrace(x_1,\cdots, x_n)\in(\R^2)^n \phantom{a}|\phantom{a} x_i\neq x_j\text{ for every }i\neq j\rbrace\] 
    
    and $\Sigma_\gamma<\Sigma_n$ is a subgroup of the Symmetric group of degree $n$ which acts both on $C_2(n)$ and $\prod_{i=1}^{n}  L\Emb^{\gamma_i}$: on the former by coordinate-permutation and on the latter by permuting its factors. It is defined as the subgroup of $\Sigma_n$ that 
    preserves the quotient partition of $\{ 1,2,\cdots, n\}$ under the equivalence relation given by $i\sim j \Longleftrightarrow  L\Emb^{\gamma_i}= L\Emb^{\gamma_j}$. 
\end{remark}

\begin{remark}
In view of this result, the fundamental group of the connected component $L\Emb^\gamma(\R,\R^3)$ of a connected sum $\gamma:=\gamma_1\#\cdots\#\gamma_n$ of prime components can be expressed in terms of the fundamental group of each of the components as well as the fundamental group of the configuration space $\mathcal C_2$ and the action of the group $\Sigma_\gamma$. 
\end{remark}

\begin{remark}\label{pullingloops}
    Note that the fundamental group of $\mathcal C_2(n)$-coincides with the braid group on $n$-strands. R. Budney   
    constructs in \cite{Budney1} a homotopy equivalence between the space of long embeddings and the free little $2$-cubes object over the space of prime knots. Moreover, he explicitly gets the following short exact sequence \cite[p.19]{Budney1}:
\begin{equation}\label{SES}
0\to \pi_1\left(\mathcal{C}_2(n)\right)\times\prod_{i=1}^{n} \pi_1\left(\tilde{K}_{\gamma_i}\right) \to \pi_1\left(\tilde{K}_{\gamma}\right) \to \Sigma_\gamma\to 0
\end{equation}

where $\tilde{K}_{\gamma_i}$ and $\tilde{K}_{\gamma}$ just denote  auxilliary spaces homotopically equivalent to $L\Emb^{\gamma_i}(\R,\R^3)$ and $L\Emb^{\gamma}(\R,\R^3)$, respectively (and, thus, with the same homotopy groups).

The action of the braid group can be described in geometric terms. It follows from the work \cite{Budney1} that elements of the group $\pi_1(\tilde{K}_{\gamma})$ can be generated by the following type of isotopies. Take two adjacent long knots $\gamma_i$ and $\gamma_j$ in the long representation of $\gamma$. We can shrink $\gamma_i$ and pull it all trough the component $\gamma_j$. Likewise, we can do the same with $\gamma_j$ and $\gamma_i$ in the reverse order. If by applying these type of ``exchanges'' a finite number of times, each component $\gamma_i$ returns to its original position, then this process describes a loop coming from the action of the braid group $\pi_1\left(\mathcal C_2(n)\right)$.

Note, also, that it is possible that some of these combinations do not give raise to the initial order of the components but can yield a loop as well. This is indeed possible if some of the components were initially isotopic and they exchanged their positions. These are the loops in $\pi_1(\tilde{K}_{\gamma})$ yielding non-trivial elements in $\Sigma_\gamma$ in the short exact sequence above.

\end{remark}

\subsection{The Legendrian setting}
Consider the natural inclusion $\Gamma: L\Leg(\R,\R^3)\to L\Emb(\R,\R^3)$ of the space of long Legendrian embeddings within the space of smooth long embeddings. 

\begin{theorem}\label{StrConnSum}
    Let $\gamma$ be the connected sum of $k \geq 1$ prime long Legendrian embeddings $\gamma_1, \ldots, \gamma_k$ satisfying the following condition: $\gamma_i$ is smoothly isotopic to $\gamma_j$ if and only if they are Legendrian isotopic. If the fundamental group homomorphism induced by the natural inclusion
    
    \[
     \pi_1\left(L\Leg^{\gamma_j}(\R,\R^3)\right) \rightarrow \pi_1\left(L\Emb^{\gamma_j}(\R,\R^3)\right)
    \]
    
    is surjective for every  component $\gamma_j$,  $(j = 1, \ldots, k)$, then the analogous 
 group homomorphism
    
    \[ \pi_1\left(L\Leg^{\gamma}(\R,\R^3)\right) \rightarrow  \pi_1\left(L\Emb^\gamma(\R,\R^3)\right)
    \]
    
    is also surjective.

\begin{remark}
    Note that the condition in the hypothesis of the Theorem is satisfied for infinitely many collections of Legendrian embeddings. For instance, if  all the components are equal; i.e. $\gamma_i=\gamma_j$ for $i\neq j$, then the condition is trivially satisfied. Analogously, if all $\gamma_i, \gamma_j$ represent different smooth isotopy classes the condition is fulfilled as well. But there are infinitely many more possibilities. 
\end{remark}
    
\end{theorem}

\begin{proof}
    In view of the homotopical description of the fundamental group $\pi_1\left(L\Emb^\gamma(\R,\R^3)\right)$ associated to the connected sum long embedding $\gamma=\gamma_1\#\cdots\#\gamma_k$, it suffices to prove that we can always realise the isotopies described in Remark \ref{pullingloops} by Legendrian isotopies. 

    This follows from some of the arguments developed in \cite[5.2]{FMP4}, but let us explain it in detail for the sake of completeness. Consider two adjacent Legendrian components $\gamma_i$ and $\gamma_{i+1}$. Note that the following process is realisable by a Legendrian isotopy. First, we arbitrarily shrink the $\gamma_i$-component until it becomes small enough. Then, we can just pull this component all over the $\gamma_j$-component. If the shrinking was fine enough, the component is then arbitrarily small and this process does not produce self-intersections. See Figure \ref{Pull} for an explicit example.

\begin{figure}[h]
	\centering
	\includegraphics[width=1\textwidth]{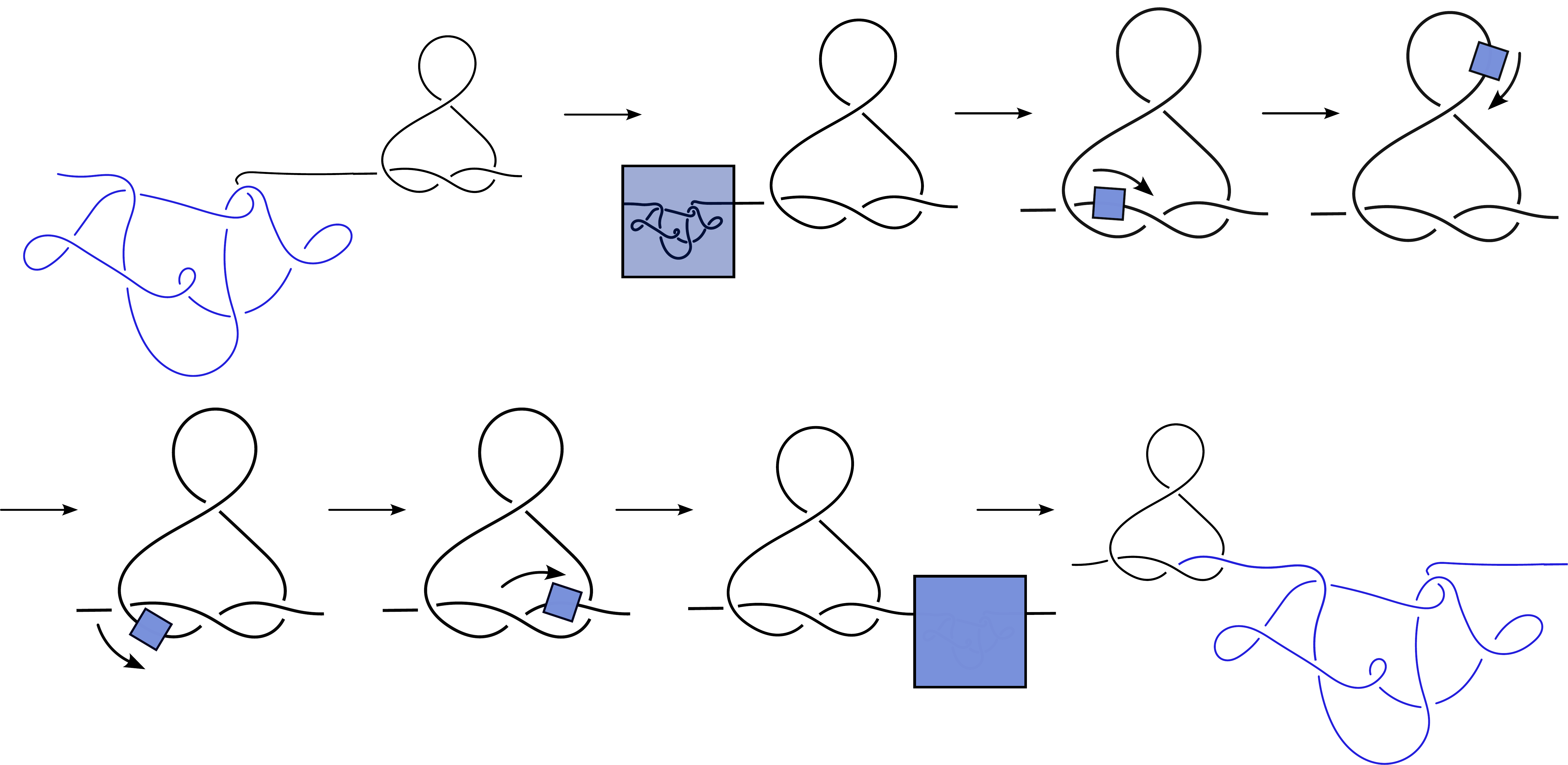}
	\caption{
 Legendrian isotopy that exchanges two adjacent prime Legendrian long embeddings: a long $m_{5_2}$-Legendrian embedding and a long Legendrian right-handed trefoil. The process is straightforward: we first shrink one of the long components, we pull it all through the other long embedding and, finally, we enlarge it until it recovers its original size.}\label{Pull}
\end{figure}

    The need of imposing the hypothesis in the statement comes from the following fact. The loop of smooth embeddings that we are trying to realise may be of the second type in Remark \ref{pullingloops}; i.e. of the ones that in the short exact sequence (\ref{SES}) have non-trivial image in $\Sigma_\gamma$.
    
    It is then clear that even if we manage to realise the aforementioned smooth isotopies by Legendrian isotopies, these may not close to a loop unless the respective components are actually Legendrian isotopic. Therefore, this is not a problem if we impose this condition from scratch, thus yielding the claim.  
\end{proof}

As a consequence of this Theorem we can finish the proof of Theorem \ref{mainPI1} by proving the statement concerning satellites.

\begin{theorem}\label{mainPISat}

The fundamental group homomorphism 
\begin{equation} \pi_1\left(\Leg(\NS^3,\xi_\std), \gamma\right)\rightarrow \pi_1\left(\Emb(\NS^1,\NS^3), \gamma\right)
\end{equation}

induced by the natural inclusion (\ref{NaturalInclusion})
is surjective for infinitely many smooth satellite embeddings and every stabilised Legendrian representative $\gamma\in\Leg(\NS^3,\xi_\std)$.
Such satellite embeddings can, moreover, be taken to be connected-sums of $k$ Legendrian hyperbolic and/or torus embeddings (for every $k\geq 2)$.
\end{theorem}

\begin{proof}
    Once again by Proposition \ref{PropSurj1}, it suffices to prove the analogous statement at the level of long-components. Take $k\in\mathbb{N}$ long embeddings satisfying the first condition from Theorem \ref{StrConnSum} and so that $
   \pi_1\left( L\Emb(\R,\R^3),\gamma\right)=\mathfrak{G}_\gamma$. Note that there exist infinitely many possibilites; e.g. $k$-copies of the same Legendrian representative of any of the embeddings in Proposition \ref{FGPI1}, Legendrian torus embeddings with the same formal invariants, etc.

   Take those $k$ long embeddings $\gamma_i$, $i=1,\cdots, k$, perform their connected-sum $\gamma=\#_{i=1}^k\gamma_i$ and add one stabilisation to $\gamma$ so that it becomes a stabilised Legendrian embedding.

\begin{figure}[h]
	\centering
	\includegraphics[width=0.6\textwidth]{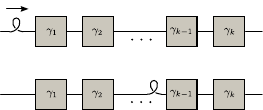}
	\caption{Schematical depiction of the connected sum $\gamma=\#_{i=1}^k\gamma_i$ of the $k$ long Legendrian embeddings considered in the proof of Theorem \ref{mainPISat}. The small kink represents the stabilisation. Note that the stabilisation can move throughout the proof in order to be recruited by each of the components (schematically depicted as grey boxes) in order to produce (through Legendrian isotopies) their respective Gramain's and Hatcher-Fox fox loops. \label{FigSum}}
\end{figure}

   The claim now just follows in the same manner as in the proof of Theorem \ref{StrConnSum} with a slight modification. Instead of assuming that each $\pi_1$-homorphism $
     \pi_1\left(L\Leg^{\gamma_j}(\R,\R^3)\right) \rightarrow \pi_1\left(L\Emb^{\gamma_j}(\R,\R^3)\right)$
     is individually surjective, we will instead realise that we can just move the stabilisation all along the embedding. This way, by Theorem \ref{Grealisations} we can realise each of the smooth fundamental groups individually (using the recruited stabilisation), one at a time. Then the claim follows by considering the additional isotopies that exchange the different components (Remark \ref{pullingloops}) as in the proof of Theorem \ref{Grealisations}.\end{proof}

\end{document}